\documentclass[12pt,leqno]{article} % Specifies the document class.
\usepackage{amsfonts}          % The preamble begins here.
\usepackage{myart98}
\title{\ \\ \ \\
\normalsize{\bf{MULTIPARAMETRIC DISSIPATIVE LINEAR STATIONARY \\
DYNAMICAL SCATTERING SYSTEMS: DISCRETE CASE, \\
II: EXISTENCE OF CONSERVATIVE
DILATIONS}}\thanks{Research supported in
part by the Ukrainian-Israeli project of scientific co-operation
(contract
no. 2M/1516-97).}}%Declares the document's title.
\author{\normalsize\bf{Dmitriy S. Kalyuzhniy}}
\date{}
\newcommand{\nspace}[2]{\ensuremath{{\mathbb{#1}}^{#2}}}
\newcommand{\Hspace}[1]{\ensuremath{\mathcal{#1}}}
\def\ifundefined#1{\expandafter\ifx\csname#1\endcsname\relax}
\newcommand{\comment}[1]{}

\usepackage{amsfonts}

\newcommand{\Cliff}[2][\comment]{\ensuremath{{\cal C}\kern-
0.18em\ell(#1,#2)}}

\ifundefined{qed}
    \DeclareMathSymbol{\qed}{0}{AMSa}{"03}
\fi
\providecommand{\eqref}[1]{\textup{(\ref{#1})}}

\begin{document}
\maketitle
\vspace{-1cm}
\begin{abstract}
\noindent
In the present paper we introduce the notion of dilation of a
multiparametric linear stationary dynamical system (systems of this
type, in particular dissipative, and conservative scattering ones were
first
introduced in \cite{K2}). We establish the criterion for existence of a
conservative dilation of a multiparametric dissipative scattering
system.
This allows to distinguish the class of so-called $N$-dissipative
systems
preserving the most important properties of one-parametric dissipative
scattering systems.
\end{abstract}
\thispagestyle{empty}
\setcounter{section}{-1}
\section{Introduction}
This paper continues the investigation of multiparametric linear
stationary dynamical systems (LSDSs), in particular dissipative,
 and conservative
scattering systems, started in \cite{K2}. Such systems represent a
generalization of LSDSs with discrete time $t\in\mathbb{Z}$, in
particular dissipative (contractive),
and conservative (unitary) scattering systems (see \cite{A} and survey
\cite{BC}) to the case $t\in\nspace{Z}{N}$. In Section~\ref{sec:pr} we
recall the necessary definitions and facts from \cite{K2}. In
Section~\ref{sec:1-dil} we recall the definition of dilation of a
one-parametric LSDS and prove the lemma in which a useful
equivalent reformulation of this definition is obtained. In
Section~\ref{sec:n-dil} the notion of dilation of a multiparametric LSDS
is
introduced, and some of its properties are discussed. In
Section~\ref{sec:cons-dil} the criterion for existence of a conservative
dilation of a multiparametric dissipative scattering LSDS is
established.
This
criterion allows to distinguish the subclass of multiparametric
dissipative
scattering LSDSs that possess conservative dilations (we call them
$N$-dissipative scattering systems) and preserve other important
properties
of one-parametric dissipative scattering LSDSs. In 
Section~\ref{sec:n-pass}
 we
characterize the class of transfer functions of all $N$-dissipative
 scattering
LSDSs with given input and output spaces as the subclass of the
generalized
Schur class of operator-valued functions on the open unit
polydisc $\nspace{D}{N}$ (the definition of this class is recalled
in Section~\ref{sec:pr}) distinguished by the condition of vanishing
at $z=0$. We prove the existence of minimal $N$-dissipative realizations
for operator-valued functions from this subclass. We establish that in
the cases $N=1$ and $N=2$ the class of $N$-dissipative scattering LSDSs
coincides with the class of all dissipative scattering LSDSs. Note that
for
$N=1$ it is a well-known result \cite{A} appearing as a system analogue
of
the classical theorem of B. Sz.-Nagy on the existence of a unitary
dilation
for an arbitrary contractive linear operator on a Hilbert space (see
\cite{Sz.-NF}). In the case $N>2$ the class of $N$-dissipative
scattering
LSDSs
turns out to be a proper subclass of the class of all dissipative
scattering
LSDSs.
\section{Preliminaries}\label{sec:pr}
In this section we shall recall some definitions and results from
\cite{K2} that will be needed in the sequel.

For $t\in
\nspace{Z}{N}$ set $|t|:=\sum
_{k=1}^Nt_k$, for each $k\in\{1,\ldots ,N\}$ let $e_k$ be the $N$-tuple
with unit on the $k$-th place and zeroes on the rest, and let
$[\Hspace{H}_1,\Hspace{H}_2]$ denote the Banach space of all linear
bounded operators mapping a separable Hilbert space $\Hspace{H}_1$ into
a separable Hilbert space $\Hspace{H}_2$. Then a \emph{multiparametric
LSDS} is, by definition, the following system of equalities:
\begin{equation}\label{eq:n-sys}
\alpha :\left\{\begin{array}{lll}
x(t)&=&\sum _{k=1}^N(A_kx(t-e_k)+B_k\phi ^-(t-e_k)),\\
\phi ^+(t)&=&\sum _{k=1}^N(C_kx(t-e_k)+D_k\phi ^-(t-e_k)),
\end{array}
\right. (|t|>0)
\end{equation}
where for all $k\in\{1,\ldots ,N\}\ A_k\in [\Hspace{X},
\Hspace{X}], B_k\in
[\Hspace{N^-},\Hspace{X}], C_k\in [\Hspace{X},\Hspace{N^+}], D_k\in
[\Hspace{N^-},\Hspace{N^+}]$, together with the initial condition
\begin{equation}\label{eq:init}
x(t)=x_0(t)
\end{equation}
where $x_0:\{t\in\nspace{Z}{N}:|t|=0\}\to\Hspace{X}$ is a
prescribed function. We call $\Hspace{X},\ \Hspace{N^-},\ \Hspace{N^+}$
respectively the \emph{state space}, the
\emph{input space} and the \emph{output space} of $\alpha$.
If one denotes the $N$-tuple of operators $T_k\quad
(k=1,\ldots ,N)$ by $\mathbf{T}:=(T_1,\ldots ,T_N)$ then for such a
system one may use the short notation $\alpha =(N; \mathbf{A}, \mathbf{B},
\mathbf{C}, \mathbf{D};\Hspace{X}, \Hspace{N^-},
\Hspace{N^+})$. Note that in the case $N=1$ a system of equalities in
\eqref{eq:n-sys} differs from the standard one (see \cite{A} or
\cite{BC})
by shift in an output signal $\phi^+$, that brings, as it is shown in
\cite{K2}, to unessential changes in the theory of one-parametric LSDSs.
The notion of dilation for this case, that will be used in the sequel,
doesn't differ from the standard one (see Section~\ref{sec:1-dil}).

Set $z\mathbf{T}:=\sum_{k=1}^Nz_kT_k$ for $N$-tuples of complex
numbers $z=(z_1,\ldots ,z_N)$ and operators $\mathbf{T}=(T_1,\ldots
,T_N)$. Then a $[\Hspace{N^-},\Hspace{N^+}]$-valued function
\begin{displaymath}
\theta_\alpha (z)=z\mathbf{D}+z\mathbf{C}{(I_\Hspace{X}-z\mathbf{A})}^{-
1}z\mathbf{B},
\end{displaymath}
which has to be considered on some neighbourhood of $z=0$ in
$\nspace{C}{N}$, is called the \emph{transfer function} of a system
$\alpha $ of the form \eqref{eq:n-sys}--\eqref{eq:init}. The system
$\alpha =(N; \mathbf{A}, \mathbf{B}, \mathbf{C}, \mathbf{D};\Hspace{X},
\Hspace{N^-}, \Hspace{N^+})$ is called
a \emph{dissipative} (resp. \emph{conservative}) \emph{scattering LSDS} if
for each $\zeta\in\nspace{T}{N}$ ($N$-dimensional torus)
\begin{displaymath}
\zeta\mathbf{G}:=\left(
\begin{array}{ll}
\zeta\mathbf{A}&\zeta\mathbf{B}\\
\zeta\mathbf{C}&\zeta\mathbf{D}
\end{array}
\right)\in [\Hspace{X}\oplus\Hspace{N^-},\Hspace{X}\oplus\Hspace{N^+}]
\end{displaymath}
is a contractive (resp. unitary) operator.
\begin{thm}\label{thm:tf-pass}
The transfer function $\theta_\alpha $ of an arbitrary dissipative
scattering LSDS $\alpha =(N; \mathbf{A},
\mathbf{B}, \mathbf{C},
\mathbf{D};\Hspace{X}, \Hspace{N^-}, \Hspace{N^+})$ belongs to the class
$B_N^0(\Hspace{N^-},\Hspace{N^+})$
consisting of all functions holomorphic on the open unit polydisc
$\nspace{D}{N}$ with contractive values from  $[\Hspace{N^-},\Hspace{N^+}]$ and
vanishing at $z=0$.
\end{thm}
Recall (see \cite{Ag}) that the \emph{generalized Schur class}
$S_N(\Hspace{N^-},\Hspace{N^+})$ is the class of functions
\begin{displaymath}
\theta (z)=\sum_{t\in\nspace{Z}{N}_+}\widehat{\theta}_tz^t
\end{displaymath}
holomorphic on $\nspace{D}{N}$ with values in $[\Hspace{N^-},\Hspace{N^+}]$
(here
$\nspace{Z}{N}_+:=\{ t\in\nspace{Z}{N}:t_k\geq 0,k=1,\ldots ,N\}$ is the
discrete positive octant, $z^t:=\prod_{k=1}^Nz_k^{t_k}$ is a usual
multipower for $t\in\nspace{Z}{N}_+$), such that for any separable
Hilbert space $\Hspace{Y}$, any $N$-tuple $\mathbf{T}=(T_1,\ldots ,T_N)$
of commuting contractions on $\Hspace{Y}$ and for any positive $r<1$ one
has
\begin{displaymath}
\| \theta (r\mathbf{T})\|\leq 1
\end{displaymath}
where
\begin{displaymath}
\theta (r\mathbf{T})=\theta (rT_1,\ldots
,rT_N) := \sum_{t\in\nspace{Z}{N}_+}\widehat{\theta _t}\otimes
{(r\mathbf{T})}^t\in [\Hspace{N^-}\otimes \Hspace{Y},\Hspace{N^+}\otimes
\Hspace{Y}]
\end{displaymath}
(the convergence of this series is understood in the
sense of norm in the Banach space
$[\Hspace{N^-}\otimes \Hspace{Y},\Hspace{N^+}\otimes \Hspace{Y}]$). If
$N=1$ then due to the von Neumann inequality (see \cite{vN}) we have
$S_N(\Hspace{N^-},\Hspace{N^+})=S(\Hspace{N^-},\Hspace{N^+})$ i.e. the
\emph{Schur class} consisting of all functions holomorphic on the 
open unit disc  $\mathbb{D}$ with contractive values from
$[\Hspace{N^-},\Hspace{N^+}]$.

Denote by $S_N^0(\Hspace{N^-},\Hspace{N^+})$ the subclass of those
functions from $S_N(\Hspace{N^-},\Hspace{N^+})$ that vanish at $z=0$.
\begin{thm}\label{thm:tf-cons}
The class of transfer functions of $N$-parametric conservative
scattering LSDSs with the input space $\Hspace{N^-}$ and the output
space $\Hspace{N^+}$ coincides with $S_N^0(\Hspace{N^-},\Hspace{N^+})$.
\end{thm}
In conclusion of this section let us remark that it is not difficult to
verify the following inclusion:
\begin{equation}\label{eq:inclus}
S_N^0(\Hspace{N^-},\Hspace{N^+})\subseteq B_N^0(\Hspace{N^-
},\Hspace{N^+}).
\end{equation}
It is known that for $N=1$ (see \cite{vN}) and for $N=2$
(see \cite{An}) we have in fact the sign ``='' in \eqref{eq:inclus} for
any
$\Hspace{N^-}$ and $\Hspace{N^+}$, i.e. the classes $S_N^0(\Hspace{N^-
},\Hspace{N^+})$ and $B_N^0(\Hspace{N^-},\Hspace{N^+})$ coincide. For
$N>2$, as it follows from \cite{V}, these classes do not coincide, i.e.
we have the strict inclusion in \eqref{eq:inclus} for
any $\Hspace{N^-}$ and $\Hspace{N^+}$ different from $\{ 0\}$.

\section{Lemma on dilations of one-parametric LSDSs}\label{sec:1-dil}

In this section the question is the case $N=1$, i.e. one-parametric
systems of the form
\begin{displaymath}
\alpha :\left\{ \begin{array}{lll}
x(t)&=&Ax(t-1)+B\phi ^-(t-1),\\
\phi ^+(t)&=&Cx(t-1)+D\phi ^-(t-1),
\end{array}
\right.\quad
(t=1,2,\ldots )
\end{displaymath}
where $A\in [\Hspace{X},
\Hspace{X}], B\in
[\Hspace{N^-},\Hspace{X}], C\in [\Hspace{X},\Hspace{N^+}], D\in
[\Hspace{N^-},\Hspace{N^+}]$, and initial condition will be unessential
for our consideration; we shall write $\alpha =(A, B, C, D; \Hspace{X},
\Hspace{N^-}, \Hspace{N^+}):=(1; \mathbf{A}, \mathbf{B}, \mathbf{C},\\
\mathbf{D};\Hspace{X}, \Hspace{N^-}, \Hspace{N^+})$. Recall (see e.g.
\cite{A} or \cite{BC}) that the LSDS $\widetilde{\alpha}
=(\widetilde{A},\widetilde{B},\widetilde{C},D;\widetilde{\Hspace{X}},
\Hspace{N^-}, \Hspace{N^+})$ is said to be a \emph{dilation of} the
\emph{LSDS} $\alpha =(A, B, C, D; \Hspace{X}, \Hspace{N^-},
\Hspace{N^+})$ if there are such subspaces $\Hspace{D}$ and
$\Hspace{D}_*$ in $\widetilde{\Hspace{X}}$ that
\begin{equation}\label{eq:ort-sum}
\widetilde{\Hspace{X}}=\Hspace{D}\oplus\Hspace{X}\oplus\Hspace{D_*},
\end{equation}
\begin{equation}\label{eq:inv}
\widetilde{A}\Hspace{D}\subset\Hspace{D},\quad
\widetilde{C}\Hspace{D}=\{ 0\} ,\quad
{\widetilde{A}}^*\Hspace{D_*}\subset\Hspace{D_*},\quad
\widetilde{B}^*\Hspace{D_*}=\{ 0\},
\end{equation}
\begin{equation}\label{eq:proj}
A=P_\Hspace{X}\widetilde{A}|\Hspace{X},\quad
B=P_\Hspace{X}\widetilde{B},\quad
C=\widetilde{C}|\Hspace{X}
\end{equation}
(here $P_\Hspace{X}$ denotes the orthogonal projector onto
$\Hspace{X}$ in $\widetilde{\Hspace{X}}$).
\begin{lem}\label{lem:1-dil}
The LSDS $\widetilde{\alpha}
=(\widetilde{A},\widetilde{B},\widetilde{C},D;\widetilde{\Hspace{X}},
\Hspace{N^-}, \Hspace{N^+})$ is a dilation of the LSDS $\alpha =(A, B,\\
C, D; \Hspace{X}, \Hspace{N^-}, \Hspace{N^+})$ if and only if
$\Hspace{X}\subset\widetilde{\Hspace{X}}$ and for all $n\in\mathbb{Z}_+$
the following equalities hold:
\begin{equation}\label{eq:mom}
A^n=P_\Hspace{X}{\widetilde{A}}^n|X,\
A^nB=P_\Hspace{X}{\widetilde{A}}^n\widetilde{B},\
CA^n=\widetilde{C}{\widetilde{A}}^n|\Hspace{X},\
CA^nB=\widetilde{C}{\widetilde{A}}^n\widetilde{B}.
\end{equation}
\end{lem}
\begin{proof}
Suppose that $\widetilde{\alpha}$ is a dilation of $\alpha$. Then
$\Hspace{X}\subset\widetilde{\Hspace{X}}$, and by \eqref{eq:proj} for
$n=1$ the first equality in \eqref{eq:mom} holds (note that for $n=0$ it
holds trivially), and for $n=0$ the second and the third equalities in
\eqref{eq:mom} hold. Thus we have the base of induction for the proof of
the first three relations in \eqref{eq:mom}. Let for
$n=k\in\mathbb{Z}_+$ these relations are true. We will show that then
for $n=k+1$ they are also true. We have
\begin{eqnarray*}
A^{k+1} & = & A\cdot
A^k=P_\Hspace{X}\widetilde{A}P_\Hspace{X}{\widetilde{A}}^k|\Hspace{X}=
P_\Hspace{X}\widetilde{A}(I_{\widetilde{\Hspace{X}}}-P_\Hspace{D}-
P_{\Hspace{D}_*}){\widetilde{A}}^k|\Hspace{X} \\
        & = & P_\Hspace{X}{\widetilde{A}}^{k+1}|\Hspace{X}-
P_\Hspace{X}\widetilde{A}P_\Hspace{D}{\widetilde{A}}^k|\Hspace{X}-
P_\Hspace{X}\widetilde{A}P_{\Hspace{D}_*}{\widetilde{A}}^k|\Hspace{X} =
P_\Hspace{X}{\widetilde{A}}^{k+1}|\Hspace{X}
\end{eqnarray*}
since, according to \eqref{eq:ort-sum} and \eqref{eq:inv},
$P_\Hspace{X}\widetilde{A}P_\Hspace{D}=0$ and
$P_{\Hspace{D}_*}{\widetilde{A}}^k|\Hspace{X}=0$;
\begin{eqnarray*}
A^{k+1}B & = & A\cdot A^kB =
P_\Hspace{X}\widetilde{A}P_\Hspace{X}{\widetilde{A}}^k\widetilde{B} =
P_\Hspace{X}\widetilde{A}(I_{\widetilde{\Hspace{X}}}-P_\Hspace{D}-
P_{\Hspace{D}_*}){\widetilde{A}}^k\widetilde{B} \\
        & = & P_\Hspace{X}{\widetilde{A}}^{k+1}\widetilde{B}-
P_\Hspace{X}\widetilde{A}P_\Hspace{D}{\widetilde{A}}^k\widetilde{B}-
P_\Hspace{X}\widetilde{A}P_{\Hspace{D}_*}{\widetilde{A}}^k\widetilde{B}
=
P_\Hspace{X}{\widetilde{A}}^{k+1}\widetilde{B}
\end{eqnarray*}
since, according to \eqref{eq:ort-sum} and \eqref{eq:inv},
$P_\Hspace{X}\widetilde{A}P_\Hspace{D}=0$ and
$P_{\Hspace{D}_*}{\widetilde{A}}^k\widetilde{B}=0$;
\begin{eqnarray*}
CA^{k+1} & = & CA^k\cdot A =
\widetilde{C}{\widetilde{A}}^kP_\Hspace{X}\widetilde{A}|\Hspace{X}=
\widetilde{C}{\widetilde{A}}^k(I_{\widetilde{\Hspace{X}}}-P_\Hspace{D}-
P_{\Hspace{D}_*})\widetilde{A}|\Hspace{X} \\
        & = & \widetilde{C}{\widetilde{A}}^{k+1}|\Hspace{X}-
\widetilde{C}{\widetilde{A}}^kP_\Hspace{D}\widetilde{A}|\Hspace{X}-
\widetilde{C}{\widetilde{A}}^k
P_{\Hspace{D}_*}\widetilde{A}|\Hspace{X} =
\widetilde{C}{\widetilde{A}}^{k+1}|\Hspace{X}
\end{eqnarray*}
since, according to \eqref{eq:ort-sum} and \eqref{eq:inv},
$\widetilde{C}{\widetilde{A}}^kP_\Hspace{D}=0$ and
$P_{\Hspace{D}_*}\widetilde{A}|\Hspace{X}=0$.
Thus we established by induction on $n$ that the first three relations
in \eqref{eq:mom} hold. We get from here for an arbitrary
$n\in\mathbb{Z}_+$
\begin{eqnarray*}
CA^nB & = & (\widetilde{C}{\widetilde{A}}^n|\Hspace{X})\cdot
(P_\Hspace{X}\widetilde{B}) =
\widetilde{C}{\widetilde{A}}^nP_\Hspace{X}\widetilde{B} =
\widetilde{C}{\widetilde{A}}^n(I_{\widetilde{\Hspace{X}}}-P_\Hspace{D}-
P_{\Hspace{D}_*})\widetilde{B} \\
        & = & \widetilde{C}{\widetilde{A}}^n\widetilde{B}-
\widetilde{C}{\widetilde{A}}^nP_\Hspace{D}\widetilde{B}-
\widetilde{C}{\widetilde{A}}^nP_{\Hspace{D}_*}\widetilde{B} =
\widetilde{C}{\widetilde{A}}^n\widetilde{B}
\end{eqnarray*}
since, according to \eqref{eq:ort-sum} and \eqref{eq:inv},
$\widetilde{C}{\widetilde{A}}^nP_\Hspace{D}=0$ and
$P_{\Hspace{D}_*}\widetilde{B}=0$. Thus the fourth relation in
\eqref{eq:mom} is also true for any $n\in\mathbb{Z}_+$.

Conversely, let $\Hspace{X}\subset\widetilde{\Hspace{X}}$ and for all
$n\in\mathbb{Z}_+$ the equalities in \eqref{eq:mom} hold. Then set
\begin{displaymath}
\Hspace{D}:=\bigvee_{n=0}^\infty{\widetilde{A}}^n((\widetilde{A}-
A)\Hspace{X}+(\widetilde{B}-B)\Hspace{N^-})
\end{displaymath}
where the symbol ``$\bigvee$'' denotes the closure of a linear span of
some lineals, $\Hspace{U}+\Hspace{V}:=\{u+v:u\in\Hspace{U},
v\in\Hspace{V}\}$ is the sum of lineals $\Hspace{U}$ and $\Hspace{V}$ in
some space, $(\widetilde{A}-
A)\Hspace{X}:=\{\widetilde{x}\in\widetilde{\Hspace{X}}:\widetilde{x}=
\widetilde{A}x-Ax, x\in\Hspace{X}\},\ (\widetilde{B}-B)\Hspace{N^-
}:=\{\widetilde{x}\in\widetilde{\Hspace{X}}:\widetilde{x}=\widetilde{B}
\phi^--B\phi^-, \phi^-\in\Hspace{N^-}\}$. Then
$\Hspace{D}\perp\Hspace{X}$. Indeed, for arbitrary $x\in\Hspace{X},\
\phi^-\in\Hspace{N^-}$ and $n\in\mathbb{Z_+}$ we have
\begin{eqnarray*}
\lefteqn{P_\Hspace{X}{\widetilde{A}}^n((\widetilde{A}-
A)x+(\widetilde{B}-
B)\phi^-)} \\
&=& P_\Hspace{X}{\widetilde{A}}^{n+1}x-
(P_\Hspace{X}{\widetilde{A}}^n|\Hspace{X})\cdot
(Ax)+P_\Hspace{X}{\widetilde{A}}^n\widetilde{B}\phi^--
(P_\Hspace{X}{\widetilde{A}}^n|\Hspace{X})\cdot (B\phi^-) \\
&=& A^{n+1}x-A^n\cdot Ax+A^nB\phi^--A^n\cdot B\phi^-=0.
\end{eqnarray*}
Hence $P_\Hspace{X}\Hspace{D}=\{ 0\}$, and $\Hspace{D}\perp\Hspace{X}$.
Set
\begin{displaymath}
\Hspace{D}_*:=\widetilde{\Hspace{X}}\ominus(\Hspace{X}\oplus\Hspace{D}).
\end{displaymath}
Then \eqref{eq:ort-sum} is valid. From the definition of $\Hspace{D}$ we
obtain that $\widetilde{A}\Hspace{D}\subset\Hspace{D}$. Further, for
arbitrary $x\in\Hspace{X},\ \phi^-\in\Hspace{N^-}$ and
$n\in\mathbb{Z_+}$ we have
\begin{eqnarray*}
\lefteqn{\widetilde{C}{\widetilde{A}}^n((\widetilde{A}-
A)x+(\widetilde{B}-
B)\phi^-)} \\
&=& \widetilde{C}{\widetilde{A}}^{n+1}x-
(\widetilde{C}{\widetilde{A}}^n|\Hspace{X})\cdot
(Ax)+\widetilde{C}{\widetilde{A}}^n\widetilde{B}\phi^--
(\widetilde{C}{\widetilde{A}}^n|\Hspace{X})\cdot (B\phi^-) \\
&=& CA^{n+1}x-CA^n\cdot Ax+CA^nB\phi^--CA^n\cdot B\phi^-=0.
\end{eqnarray*}
>From here we obtain that $\widetilde{C}\Hspace{D}=\{ 0\}$.
For an arbitrary $x\in\Hspace{X}$ we have
\begin{displaymath}
\widetilde{A}x=(\widetilde{A}x-Ax)+Ax\in\overline{(\widetilde{A}-
A)\Hspace{X}}\oplus\Hspace{X}\subset\Hspace{D}\oplus\Hspace{X}
\end{displaymath}
(here $\overline{\Hspace{U}}$ denotes the closure of $\Hspace{U}$). It was
shown above
that $\widetilde{A}\Hspace{D}\subset\Hspace{D}$. Hence
$\widetilde{A}(\Hspace{D}\oplus\Hspace{X})\subset\Hspace{D}\oplus
\Hspace{X}$. From here we get
${\widetilde{A}}^*\Hspace{D}_*={\widetilde{A}}^*{(\Hspace{D}\oplus
\Hspace{X})}^\perp\subset{(\Hspace{D}\oplus\Hspace{X})}^\perp =\Hspace{D}_*$.
For
an arbitrary $\phi^-\in\Hspace{N^-}$ we have
\begin{displaymath}
\widetilde{B}\phi^-=(\widetilde{B}\phi^--B\phi^-)+B\phi^-
\in\overline{(\widetilde{B}-B)\Hspace{N^-
}}\oplus\Hspace{X}\subset\Hspace{D}\oplus\Hspace{X}={(\Hspace{D}_*)}^
\perp.
\end{displaymath}
>From here we get $\widetilde{B}^*\Hspace{D}_*=\{ 0\}$. Thus relations in
\eqref{eq:inv} are true. The equalities in \eqref{eq:proj} are the
special cases of the equalities in\eqref{eq:mom}. Finally, we have
obtained that $\widetilde{\alpha }
=(\widetilde{A},\widetilde{B},\widetilde{C}, D;\widetilde{\Hspace{X}},
\Hspace{N^-},\\ \Hspace{N^+})$ is a dilation of $\alpha =(A, B, C, D;
\Hspace{X}, \Hspace{N^-}, \Hspace{N^+})$.
\end{proof}
\begin{rem}\label{rem:dil-op}
>From Lemma~\ref{lem:1-dil}, in particular, the well-known result (see
e.g. \cite{A}) follows: if $\widetilde{\alpha }
=(\widetilde{A},\widetilde{B},\widetilde{C}, D; \widetilde{\Hspace{X}},
\Hspace{N^-}, \Hspace{N^+})$ is a dilation of $\alpha =(A, B, C, D;
\Hspace{X}, \Hspace{N^-}, \Hspace{N^+})$ then
\begin{displaymath}
\forall n\in\mathbb{Z_+}\quad
A^n=P_\Hspace{X}{\widetilde{A}}^n|\Hspace{X},
\end{displaymath}
i.e. the \emph{main operator} $\widetilde{A}$ \emph{of} the
\emph{system} $\widetilde{\alpha }$ is a \emph{dilation of} the main
\emph{operator} $A$ of  the system $\alpha $.
\end{rem}

\section{The notion of dilation of a multiparametric LSDS}
\label{sec:n-dil}
\begin{defn}\label{defn:n-dil}
We shall call the LSDS $\widetilde{\alpha} =(N; \widetilde{\mathbf{A}},
\widetilde{\mathbf{B}}, \widetilde{\mathbf{C}},
\mathbf{D};\widetilde{\Hspace{X}}, \Hspace{N^-}, \Hspace{N^+})$ a
\emph{dilation of} the \emph{multiparametric LSDS} $\alpha =(N;
\mathbf{A}, \mathbf{B}, \mathbf{C}, \mathbf{D};\Hspace{X}, \Hspace{N^-},
\Hspace{N^+})$ if for each $\zeta\in\nspace{T}{N}$ the one-parametric
LSDS $\widetilde{\alpha}_\zeta
:=(\zeta\widetilde{\mathbf{A}},\zeta\widetilde{\mathbf{B}},
\zeta\widetilde{\mathbf{C}},\zeta\mathbf{D};\widetilde{\Hspace{X}},
\Hspace{N^-}, \Hspace{N^+})$ is a dilation of the one-parametric LSDS
$\alpha_\zeta :
=(\zeta\mathbf{A},\zeta\mathbf{B},\zeta\mathbf{C},\zeta\mathbf{D};
\Hspace{X}, \Hspace{N^-}, \Hspace{N^+})$, i.e. for each
$\zeta\in\nspace{T}{N}$ there are such subspaces $\Hspace{D}_\zeta$ and
$\Hspace{D}_{*,\zeta}$ in $\widetilde{\Hspace{X}}$ that
\begin{equation}\label{eq:n-ort-sum}
\widetilde{\Hspace{X}}=\Hspace{D}_\zeta\oplus\Hspace{X}\oplus
\Hspace{D}_{*,\zeta},
\end{equation}
\begin{equation}\label{eq:n-inv}
\zeta\widetilde{\mathbf{A}}\Hspace{D}_\zeta\subset\Hspace{D}_\zeta ,\
\zeta\widetilde{\mathbf{C}}\Hspace{D}_\zeta =\{ 0\},\
(\zeta\widetilde{\mathbf{A}})^*\Hspace{D}_{*,\zeta}\subset
\Hspace{D}_{*,\zeta},\
(\zeta\widetilde{\mathbf{B}})^*\Hspace{D}_{*,\zeta}=\{ 0\},
\end{equation}
\begin{equation}\label{eq:n-proj}
\zeta\mathbf{A}=P_\Hspace{X}(\zeta\widetilde{\mathbf{A}})|\Hspace{X},\
\zeta\mathbf{B}=P_\Hspace{X}(\zeta\widetilde{\mathbf{B}}),\
\zeta\mathbf{C}=(\zeta\widetilde{\mathbf{C}})|\Hspace{X}.
\end{equation}
\end{defn}
>From Lemma~\ref{lem:1-dil} we obtain the following equivalent
reformulation of Definition~\ref{defn:n-dil}.
\begin{prop}\label{prop:reform-n-dil}
The LSDS $\widetilde{\alpha} =(N; \widetilde{\mathbf{A}},
\widetilde{\mathbf{B}},\widetilde{\mathbf{C}},
\mathbf{D};\widetilde{\Hspace{X}}, \Hspace{N^-}, \Hspace{N^+})$
is a dilation of the LSDS $\alpha =(N; \mathbf{A}, \mathbf{B},
\mathbf{C}, \mathbf{D};\Hspace{X}, \Hspace{N^-}, \Hspace{N^+})$ if and
only if $\Hspace{X}\subset\widetilde{\Hspace{X}}$ and for all
$\zeta\in\nspace{T}{N}$ and $n\in\mathbb{Z_+}$ the following equalities
hold:
\begin{equation}\label{eq:n-mom}
\begin{array}{rr}
(\zeta\mathbf{A})^n=P_\Hspace{X}(\zeta\widetilde{\mathbf{A}})^n|
\Hspace{X}, & (\zeta\mathbf{A})^n\zeta\mathbf{B}=P_\Hspace{X}(\zeta
\widetilde{\mathbf{A}})^n\zeta\widetilde{\mathbf{B}},\\
\zeta\mathbf{C}(\zeta\mathbf{A})^n=\zeta\widetilde{\mathbf{C}}(\zeta
\widetilde{\mathbf{A}})^n|\Hspace{X}, &
\zeta\mathbf{C}(\zeta\mathbf{A})^n\zeta\mathbf{B}=\zeta
\widetilde{\mathbf{C}}(\zeta\widetilde{\mathbf{A}})^n\zeta
\widetilde{\mathbf{B}}.
\end{array}
\end{equation}
\end{prop}
Equating coefficients of trigonometric polynomials in $N$ variables in
both sides of equalities \eqref{eq:n-mom} we will obtain another
equivalent reformulation of Definition~\ref{defn:n-dil}, that is a
multiparametric analogue of Lemma~\ref{lem:1-dil}. For convenience of
writing of corresponding relations let us recall the notations from
\cite{K2} for the so-called \emph{symmetrized multipowers of} the
$N$-\emph{tuple} $\mathbf{A}$ and the \emph{symmetrized multipowers of}
the
$N$-\emph{tuple} $\mathbf{A}$ \emph{bordered from one side (from two
sides) by} the $N$-\emph{tuples} $\mathbf{B}$ \emph{and} $\mathbf{C}$.
If
\begin{displaymath}
c_s:=\frac{|s|!}{s_1!\cdots s_N!}\quad
(s\in\nspace{Z}{N}_+)
\end{displaymath}
denote the \emph{numbers of permutations} of $|s|$ elements of $N$
different types \emph{with repetitions} (the \emph{polynomial
coefficients)} where an element of the $j$-th type repeats itself $s_j$
times,
$[k]\ (\in\{1,\ldots ,N\})$ denotes the type of an element $k$, and
$\sigma$ runs through the set of all such permutations with repetitions,
then we set
\begin{eqnarray}
\mathbf{A}^s &:=&
c_s^{-1}\sum _\sigma A_{[\sigma
(1)]}\cdots
A_{[\sigma (|s|)]},\quad (s\in \nspace{Z}{N}_+) \label{eq:s-power} \\
(\mathbf{A\sharp
B})^s& :=& c_s^{-1}\sum _\sigma A_{[\sigma
(1)]}\cdots A_{[\sigma (|s|-1)]}B_{[\sigma (|s|)]},\quad (s\in
\nspace{Z}{N}_+\setminus
\{ 0\} )\label{eq:s-power-b} \\
(\mathbf{C\flat A})^s&:=& c_s^{-1}\sum
_\sigma C_{[\sigma
(1)]}A_{[\sigma (2)]}\cdots A_{[\sigma (|s|)]},\quad (s\in
\nspace{Z}{N}_+\setminus
\{ 0\} )\label{eq:s-power-c}\\
(\mathbf{C\flat A\sharp B})^s& :=&c_s^{-1}\sum _\sigma C_{[\sigma
(1)]}A_{[\sigma (2)]}\cdots A_{[\sigma
(|s|-1)]}B_{[\sigma (|s|)]}.\label{eq:s-power-bc} \\  
&& (s\in \nspace{Z}{N}_+\setminus \{ 0,
e_1,\ldots, e_N\} ) \nonumber
\end{eqnarray}
\begin{rem}\label{rem:multipower}
In case of the commutative $N$-tuple $\mathbf{A}$ we have
\begin{displaymath}
\mathbf{A}^s=\prod_{k=1}^N{A_k}^{s_k}
\end{displaymath}
i.e. a usual multipower.
\end{rem}
\begin{prop}\label{prop:power-n-dil}
The LSDS $\widetilde{\alpha} =(N; \widetilde{\mathbf{A}},
\widetilde{\mathbf{B}}, \widetilde{\mathbf{C}},
\mathbf{D};\widetilde{\Hspace{X}}, \Hspace{N^-}, \Hspace{N^+})$ is a
dilation of the LSDS $\alpha =(N; \mathbf{A}, \mathbf{B}, \mathbf{C},
\mathbf{D}; \Hspace{X}, \Hspace{N^-}, \Hspace{N^+})$ if and only if
$\Hspace{X}\subset\widetilde{\Hspace{X}}$ and the following equalities
hold:
\begin{equation}\label{eq:sym-mom}
\begin{array}{cc}
\forall s\in \nspace{Z}{N}_+ & \mathbf{A}^s =
P_\Hspace{X}\widetilde{\mathbf{A}}^s|\Hspace{X}, \\
\forall s\in\nspace{Z}{N}_+\setminus\{ 0\} & (\mathbf{A\sharp
B})^s = P_\Hspace{X}(\mathbf{\widetilde{A}\sharp
\widetilde{B}})^s,\\
\forall s\in\nspace{Z}{N}_+\setminus\{ 0\} & (\mathbf{C\flat A})^s =
(\mathbf{\widetilde{C}\flat \widetilde{A}})^s|\Hspace{X},\\
\forall s\in \nspace{Z}{N}_+\setminus \{ 0,e_1,\ldots, e_N\} &
(\mathbf{C\flat A\sharp B})^s = (\mathbf{\widetilde{C}\flat
\widetilde{A}\sharp \widetilde{B}})^s.
\end{array}
\end{equation}
\end{prop}
\begin{rem}
The equalities in the first line of \eqref{eq:sym-mom} mean that the
$N$-\emph{tuple}
$\widetilde{\mathbf{A}}$ \emph{of main operators of} the \emph{system}
$\widetilde{\alpha}$ is, by definition, a \emph{dilation of} the 
$N$-\emph{tuple} $\mathbf{A}$ of main operators of the system $\alpha $ (cf.
Remark~\ref{rem:dil-op}). In case of the commutative $N$-tuples
$\mathbf{A}$ and $\widetilde{\mathbf{A}}$ this coincides with the
definition of dilation for $N$-tuples of operators by
\cite{Sz.-NF} (see Remark~\ref{rem:multipower}).
\end{rem}
\begin{rem}\label{rem:z-n-dil}
It follows from Proposition~\ref{prop:power-n-dil} that one can replace
$\zeta\in\nspace{T}{N}$ by $z\in\nspace{C}{N}$ in \eqref{eq:n-mom} from
Proposition~\ref{prop:reform-n-dil}, and in
\eqref{eq:n-ort-sum}--\eqref{eq:n-proj} from
Definition~\ref{defn:n-dil}, thus the LSDS
$\widetilde{\alpha }$ is a dilation of the LSDS $\alpha $ if and only if
for each $z\in\nspace{C}{N}$ there are such subspaces $\Hspace{D}_z$ and
$\Hspace{D}_{*,z}$ in $\widetilde{\Hspace{X}}$ that
\begin{equation}\label{eq:n-ort-sum'}
\widetilde{\Hspace{X}}=\Hspace{D}_z\oplus\Hspace{X}\oplus
\Hspace{D}_{*,z},
\end{equation}
\begin{equation}\label{eq:n-inv'}
z\widetilde{\mathbf{A}}\Hspace{D}_z\subset\Hspace{D}_z ,\
z\widetilde{\mathbf{C}}\Hspace{D}_z =\{ 0\},\
(z\widetilde{\mathbf{A}})^*\Hspace{D}_{*,z}\subset
\Hspace{D}_{*,z},\
(z\widetilde{\mathbf{B}})^*\Hspace{D}_{*,z}=\{ 0\},
\end{equation}
\begin{equation}\label{eq:n-proj'}
z\mathbf{A}=P_\Hspace{X}(z\widetilde{\mathbf{A}})|\Hspace{X},\
z\mathbf{B}=P_\Hspace{X}(z\widetilde{\mathbf{B}}),\
z\mathbf{C}=(z\widetilde{\mathbf{C}})|\Hspace{X}.
\end{equation}
\end{rem}
\begin{rem}
Symmetrized multipowers that were defined in
\eqref{eq:s-power}--\eqref{eq:s-power-bc} take part in expressions for
states $x(t)$
and output signals $\phi^+(t)$ of a multiparametric LSDS $\alpha $
through states $x_0(\tau)$ from \eqref{eq:init} and input signals
$\phi^-(\tau)$ at \emph{preceding} to $t$ \emph{moments} $\tau\neq t$ of
``multidimensional time'' (we set $\tau\leq t$ if $t-
\tau\in\nspace{Z}{N}_+$), that are deduced from the recurrent relations
from \eqref{eq:n-sys} and the initial condition \eqref{eq:init} (see
\cite{K2}). Thus the algebraic definition of dilation from
Proposition~\ref{prop:power-n-dil} is connected with consideration of
system in ``multidimensional time'' domain, whereas the initial
geometric Definition~\ref{defn:n-dil} is connected with considerations
in ``multidimensional frequency'' domain or with so-called $Z$-transform
$\widehat{\alpha }$ of a system $\alpha $ (see Remark~\ref{rem:z-n-dil}
and \cite{K2}).
\end{rem}
\begin{prop}\label{prop:tf-n-dil}
The transfer functions of the system $\alpha =(N; \mathbf{A},
\mathbf{B}, \mathbf{C}, \mathbf{D};
\Hspace{X}, \Hspace{N^-}, \Hspace{N^+})$ and of its dilation
$\widetilde{\alpha} =(N;
\widetilde{\mathbf{A}},\widetilde{\mathbf{B}},\widetilde{\mathbf{C}},
\mathbf{D};\widetilde{\Hspace{X}}, \Hspace{N^-}, \Hspace{N^+})$
coincide.
\end{prop}
\begin{proof}
The transfer functions of $\alpha $ and $\widetilde{\alpha }$
\begin{equation}\label{eq:tf}
\theta_\alpha (z)=z\mathbf{D}+z\mathbf{C}(I_{\Hspace{X}}-z\mathbf{A})^{-
1}z\mathbf{B}=z\mathbf{D}+\sum_{n=0}^\infty
z\mathbf{C}(z\mathbf{A})^nz\mathbf{B},
\end{equation}
\begin{equation}\label{eq:tf-tilda}
\theta_{\widetilde{\alpha}}
(z)=z\mathbf{D}+z\widetilde{\mathbf{C}}(I_{\widetilde{\Hspace{X}}}-
z\widetilde{\mathbf{A}})^{-
1}z\widetilde{\mathbf{B}}=z\mathbf{D}+\sum_{n=0}^\infty
z\widetilde{\mathbf{C}}(z\widetilde{\mathbf{A}})^nz
\widetilde{\mathbf{B}}
\end{equation}
are defined and holomorphic on some neighbourhood of $z=0$ in
$\nspace{C}{N}$. In particular, the series in \eqref{eq:tf-tilda}
converges to $\theta_{\widetilde{\alpha }} (z)$ in operator norm
uniformly and absolutely on compact subsets of the domain
$\Upsilon :=\{ z\in\nspace{C}{N}:\|z\widetilde{\mathbf{A}}\|<1\}$. In
this domain $\|z\mathbf{A}\| =\|
P_\Hspace{X}(z\widetilde{\mathbf{A}})|\Hspace{X}\|
\leq\|z\widetilde{\mathbf{A}}\| <1$. Therefore the series in
\eqref{eq:tf} converges to $\theta_\alpha (z)$ in operator norm
uniformly and absolutely on compact subsets of $\Upsilon$. Besides, it
follows from Proposition~\ref{prop:reform-n-dil} and
Remark~\ref{rem:z-n-dil} that for all $z\in \nspace{C}{N}$ and
$n\in\mathbb{Z_+}\quad
z\mathbf{C}(z\mathbf{A})^nz\mathbf{B}=z\widetilde{\mathbf{C}}(z
\widetilde{\mathbf{A}})^nz\widetilde{\mathbf{B}}$, and hence for all
$z\in\Upsilon$, according to \eqref{eq:tf} and \eqref{eq:tf-tilda}, we
have $\theta_\alpha (z)=\theta_{\widetilde{\alpha}} (z)$.
\end{proof}
\begin{defn}
We shall call a multiparametric LSDS \emph{minimal} if it is not a
dilation of any system other than itself.
\end{defn}
\begin{prop}\label{prop:exist-min}
For an arbitrary LSDS $\alpha =(N; \mathbf{A}, \mathbf{B}, \mathbf{C},
\mathbf{D};\Hspace{X}, \Hspace{N^-}, \Hspace{N^+})$ there exists a minimal LSDS
$\alpha_{min} =(N; \mathbf{A}_{min}, \mathbf{B}_{min},
\mathbf{C}_{min}, \mathbf{D};\Hspace{X}_{min}, \Hspace{N^-},
\Hspace{N^+})$ such that $\alpha $ is a dilation of $\alpha _{min}$.
\end{prop}
\begin{proof}
We will use the Zorn lemma (see e.g. \cite{M}). Consider the set
$\Sigma_\alpha $ of all systems $\alpha _\gamma,\ \gamma\in\Gamma$ (here
$\Gamma$ is some set of indices), for which $\alpha $ is a dilation.
Then $\Sigma_\alpha $ is a partially ordered set with respect to the
relation ``$\succ$'': we shall write $\alpha _{\gamma_1}\succ\alpha
_{\gamma_2}$ if $\alpha _{\gamma_1}$ is a dilation of
$\alpha _{\gamma_2}$. For the existence in $\Sigma_\alpha$ of a minimal
element (which is a minimal system with the dilation $\alpha $) it is
sufficient to prove that any chain $\mathfrak{C}_\alpha $ in
$\Sigma_\alpha $ has a lower bound. Without loss of generality one can
suppose that $\mathfrak{C}_\alpha $ contains the element
$\alpha _0=\alpha $:
\begin{displaymath}
\mathfrak{C}_\alpha :\quad \alpha =\alpha _0\succ\ldots \succ
\alpha _\gamma\ldots
\end{displaymath}
(the directed set of indices for this chain will be denoted by $\Gamma_0$). If
$\Gamma_0$ is finite then $\mathfrak{C}_\alpha $ has the 
minimal element $\alpha _{\gamma_*}$ which is a desired lower bound for
$\mathfrak{C}_\alpha $ in this case. Now let the directed set $\gamma_0$ 
be infinite. Evidently, the corresponding state spaces for systems from
$\mathfrak{C}_\alpha $ are completely ordered by inclusion ``$\supseteq$'', i.e.
we obtain the
chain
\begin{displaymath}
\mathfrak{C}_\Hspace{X}:\quad \Hspace{X}=\Hspace{X}_0\supseteq\ldots
\supseteq\Hspace{X}_\gamma\supseteq\ldots .
\end{displaymath}
Set $\Hspace{X}_*:=\bigcap_{\gamma\in\Gamma_0}\Hspace{X}_\gamma$. Then
(see e.g. \cite{AkG})
\begin{displaymath}
P_{\Hspace{X}_*}=\mbox{s-}\lim_{\gamma\in\Gamma_0}P_{\Hspace{X}_\gamma}.
\end{displaymath}
Set $\alpha_*:=(N; \mathbf{A}_*, \mathbf{B}_*, \mathbf{C}_*,
\mathbf{D};\Hspace{X}_*, \Hspace{N^-}, \Hspace{N^+})$ where the
$N$-tuples $\mathbf{A}_*, \mathbf{B}_*, \mathbf{C}_*$ of operators are
defined by formulas:
\begin{displaymath}
A_{*,k}=P_{\Hspace{X}_*}A_k|\Hspace{X}_*,\ B_{*,k}=P_{\Hspace{X}_*}B_k,\
C_{*,k}=C_k|\Hspace{X}_*.\ (k=1,\ldots ,N)
\end{displaymath}
Then $\alpha $ is a dilation of $\alpha _*$, i.e.
$\alpha _*\in\Sigma_\alpha $. To show this we shall verify the
equalities in \eqref{eq:n-mom} for these two systems and imply
Proposition~\ref{prop:reform-n-dil}. For arbitrary
$\zeta\in\nspace{T}{N}$ and $n\in\mathbb{Z_+}$ we have
\begin{eqnarray*}
(\zeta\mathbf{A}_*)^n &=&
(P_{\Hspace{X}_*}(\zeta\mathbf{A}))^n|\Hspace{X}_*=
\mbox{s-
}\lim_{\gamma\in\Gamma_0}(P_{\Hspace{X}_\gamma}(\zeta\mathbf{A}))^n|
\Hspace{X}_* \\
&=& \mbox{s-
}\lim_{\gamma\in\Gamma_0}(P_{\Hspace{X}_\gamma}(\zeta\mathbf{A})|
\Hspace{X}_\gamma)^n|\Hspace{X}_*
= \mbox{s-}\lim_{\gamma\in\Gamma_0}(\zeta\mathbf{A}_\gamma)^n|
\Hspace{X}_* \\
&=& \mbox{s-
}\lim_{\gamma\in\Gamma_0}P_{\Hspace{X}_\gamma}(\zeta\mathbf{A})^n|
\Hspace{X}_*=
P_{\Hspace{X}_*}(\zeta\mathbf{A})^n|
\Hspace{X}_*,\\
(\zeta\mathbf{A}_*)^n\zeta\mathbf{B}_* &=&
(P_{\Hspace{X}_*}(\zeta\mathbf{A}))^nP_{\Hspace{X}_*}(\zeta
\mathbf{B})=
\mbox{s-
}\lim_{\gamma\in\Gamma_0}(P_{\Hspace{X}_\gamma}(\zeta\mathbf{A}))^n
P_{\Hspace{X}_\gamma}(\zeta\mathbf{B}) \\
&=& \mbox{s-
}\lim_{\gamma\in\Gamma_0}(\zeta\mathbf{A}_\gamma)^n\zeta\mathbf{B}_
\gamma =\mbox{s-}\lim_{\gamma\in\Gamma_0}P_{\Hspace{X}_
\gamma}(\zeta\mathbf{A})^n\zeta\mathbf{B}=P_{\Hspace{X}_*}(\zeta
\mathbf{A})^n\zeta\mathbf{B},
\end{eqnarray*}
and other equalities in \eqref{eq:n-mom} are verified analogously.

Now let us show that for each $\gamma\in\Gamma_0\quad\alpha _\gamma$ is
a dilation of $\alpha _*$. Indeed, for any $\gamma\in\Gamma_0$ and
$n\in\mathbb{Z_+}$ we have
\begin{eqnarray*}
(\zeta\mathbf{A}_*)^n=P_{\Hspace{X}_*}(\zeta\mathbf{A})^n|\Hspace{X}_*=
P_{\Hspace{X}_*}(P_{\Hspace{X}_\gamma}(\zeta\mathbf{A})^n|\Hspace{X}_
\gamma)|\Hspace{X}_*=P_{\Hspace{X}_*}(\zeta\mathbf{A}_\gamma)^n|
\Hspace{X}_* ,\\
(\zeta\mathbf{A}_*)^n\zeta\mathbf{B}_*=P_{\Hspace{X}_*}(\zeta\mathbf{A})
^n\zeta\mathbf{B}=
P_{\Hspace{X}_*}P_{\Hspace{X}_\gamma}(\zeta\mathbf{A})^n\zeta\mathbf{B}
=P_{\Hspace{X}_*}(\zeta\mathbf{A}_\gamma)^n\zeta\mathbf{B}_\gamma ,\\
\zeta\mathbf{C}_*(\zeta\mathbf{A}_*)^n=\zeta\mathbf{C}(\zeta\mathbf{A})^
n|\Hspace{X}_*=
(\zeta\mathbf{C}(\zeta\mathbf{A})^n|\Hspace{X}_
\gamma)|\Hspace{X}_*=\zeta\mathbf{C}_\gamma(\zeta\mathbf{A}_\gamma)^n|
\Hspace{X}_* ,\\
\zeta\mathbf{C}_*(\zeta\mathbf{A}_*)^n\zeta\mathbf{B}_*=\zeta\mathbf{C}
(\zeta\mathbf{A})^n\zeta\mathbf{B}=
\zeta\mathbf{C}_\gamma(\zeta\mathbf{A}_\gamma)^n\zeta\mathbf{B}_\gamma .
\end{eqnarray*}
Thus $\alpha _*$ is a lower bound for $\mathfrak{C}_\alpha $, and the
proof is complete.
\end{proof}

\section{Criterion for existence of a conservative dilation of a
multiparametric dissipative scattering LSDS}\label{sec:cons-dil}
\begin{defn}
We shall say that 
$\widetilde{\alpha}
=(N; \widetilde{\mathbf{A}}, \widetilde{\mathbf{B}},
\widetilde{\mathbf{C}}, \mathbf{D}; \widetilde{\Hspace{X}},
\Hspace{N^-},
\Hspace{N^+})$ is a \emph{conservative dilation of the dissipative scattering
LSDS} $\alpha =(N;
\mathbf{A}, \mathbf{B}, \mathbf{C}, \mathbf{D};\Hspace{X}, \Hspace{N^-},
\Hspace{N^+})$ if $\widetilde{\alpha }$ is a dilation of $\alpha $,
and
$\widetilde{\alpha }$ is a conservative scattering LSDS.
\end{defn}
\begin{thm}\label{thm:criterion}
The dissipative scattering LSDS $\alpha =(N; \mathbf{A}, \mathbf{B},
\mathbf{C}, \mathbf{D};\Hspace{X}, \Hspace{N^-}, \Hspace{N^+})$ allows
a conservative dilation if and only if the corresponding linear
operator-valued function
\begin{equation}\label{eq:linear-func}
L_\mathbf{G}(z):=z\mathbf{G}=\left(
\begin{array}{ll}
z\mathbf{A} & z\mathbf{B} \\
z\mathbf{C} & z\mathbf{D}
\end{array}
\right)\quad (z\in\nspace{D}{N})
\end{equation}
belongs to the class $S_N^0(\Hspace{X}\oplus
\Hspace{N^-},\Hspace{X}\oplus\Hspace{N^+})$ (the definition of this
class can be found in Section~\ref{sec:pr}).
\end{thm}
\begin{proof}
Let the dissipative scattering LSDS $\alpha =(N; \mathbf{A}, \mathbf{B},
\mathbf{C}, \mathbf{D};\Hspace{X}, \Hspace{N^-}, \Hspace{N^+})$ possess
the conservative dilation $\widetilde{\alpha} =(N;
\widetilde{\mathbf{A}}, \widetilde{\mathbf{B}}, \widetilde{\mathbf{C}},
\mathbf{D};\widetilde{\Hspace{X}}, \Hspace{N^-}, \Hspace{N^+})$. Then
for each $\zeta\in\nspace{T}{N}$
\begin{displaymath}
\zeta\widetilde{\mathbf{G}}:=\left(
\begin{array}{ll}
\zeta\widetilde{\mathbf{A}} & \zeta\widetilde{\mathbf{B}} \\
\zeta\widetilde{\mathbf{C}} & \zeta\mathbf{D}
\end{array}
\right)\in [\widetilde{\Hspace{X}}\oplus
\Hspace{N^-},\widetilde{\Hspace{X}}\oplus\Hspace{N^+}]
\end{displaymath}
is a unitary operator. This operator allows another block partition:
\begin{displaymath}
\zeta\widetilde{\mathbf{G}}=\left(
\begin{array}{ll}
\zeta\mathbf{T} & \zeta\mathbf{F} \\
\zeta\mathbf{H} & \zeta\mathbf{S}
\end{array}
\right)\in [(\widetilde{\Hspace{X}}\ominus\Hspace{X})\oplus
(\Hspace{X}\oplus
\Hspace{N^-}),(\widetilde{\Hspace{X}}\ominus\Hspace{X})\oplus
(\Hspace{X}\oplus\Hspace{N^+})]
\end{displaymath}
where
\begin{displaymath}
\begin{array}{ll}
\zeta\mathbf{T}  =P_{\widetilde{\Hspace{X}}\ominus\Hspace{X}}(\zeta
\widetilde{\mathbf{G}})|\widetilde{\Hspace{X}}\ominus\Hspace{X}, &
\zeta\mathbf{F}=P_{\widetilde{\Hspace{X}}\ominus\Hspace{X}}(\zeta
\widetilde{\mathbf{G}})|\Hspace{X}\oplus\Hspace{N^-}, \\
\zeta\mathbf{H} =P_{\Hspace{X}\oplus\Hspace{N^+}}(\zeta
\widetilde{\mathbf{G}})|\widetilde{\Hspace{X}}\ominus\Hspace{X}, &
\zeta\mathbf{S}=P_{\Hspace{X}\oplus\Hspace{N^+}}(\zeta
\widetilde{\mathbf{G}})|\Hspace{X}\oplus\Hspace{N^-}=\zeta\mathbf{G}.
\end{array}
\end{displaymath}
It is clear that one can correspond to this partition of
$\zeta\widetilde{\mathbf{G}}$ the conservative scattering LSDS
$\beta=(N; \mathbf{T}, \mathbf{F}, \mathbf{H}, \mathbf{S};
\widetilde{\Hspace{X}}\ominus\Hspace{X}, \Hspace{X}\oplus\Hspace{N^-},
\Hspace{X}\oplus\Hspace{N^+})$ where $\mathbf{S}=\mathbf{G}$. Its
transfer function
\begin{displaymath}
\theta_\beta
(z)=z\mathbf{S}+z\mathbf{H}(I_{\widetilde{\Hspace{X}}\ominus\Hspace{X}}-
z\mathbf{T})^{-1}z\mathbf{F}
\end{displaymath}
by Theorem~\ref{thm:tf-cons} belongs to the class
$S_N^0(\Hspace{X}\oplus
\Hspace{N^-},\Hspace{X}\oplus\Hspace{N^+})$. Let us show that
$\theta_\beta (z)=z\mathbf{S}=z\mathbf{G}=L_\mathbf{G}(z)$. This will
mean that $L_\mathbf{G}\in S_N^0(\Hspace{X}\oplus
\Hspace{N^-},\Hspace{X}\oplus\Hspace{N^+})$ i.e. the necessary condition
of the present theorem. Evidently, it is sufficient to show that for any
$z\in\nspace{C}{N}$ and $n\in\mathbb{Z_+}$
\begin{displaymath}
z\mathbf{H}(z\mathbf{T})^nz\mathbf{F}=0.
\end{displaymath}
According to \eqref{eq:n-ort-sum'} we have
\begin{eqnarray*}
\lefteqn{z\mathbf{H}(z\mathbf{T})^nz\mathbf{F}} \\
&& =(P_{\Hspace{X}\oplus
\Hspace{N^+}}(z\widetilde{\mathbf{G}})|\widetilde{\Hspace{X}}\ominus
\Hspace{X})(P_{\widetilde{\Hspace{X}}\ominus\Hspace{X}}(z
\widetilde{\mathbf{G}})|\widetilde{\Hspace{X}}\ominus
\Hspace{X})^nP_{\widetilde{\Hspace{X}}\ominus\Hspace{X}}(z
\widetilde{\mathbf{G}})|\Hspace{X}\oplus\Hspace{N^-} \\
&& =\left[ \begin{array}{c}
P_\Hspace{X}(z\widetilde{\mathbf{A}}) \\
z\widetilde{\mathbf{C}}
\end{array}
\right] (P_{\Hspace{D}_z}+P_{\Hspace{D}_{*,z}})
(z\widetilde{\mathbf{A}}(P_{\Hspace{D}_z}+P_{\Hspace{D}_{*,z}}))^n
[\begin{array}{ll}
(z\widetilde{\mathbf{A}})|\Hspace{X} & z\widetilde{\mathbf{B}}
\end{array}] \\
&& =\left[ \begin{array}{c}
P_\Hspace{X}(z\widetilde{\mathbf{A}}) \\
z\widetilde{\mathbf{C}}
\end{array}
\right] P_{\Hspace{D}_{*,z}}
(z\widetilde{\mathbf{A}}(P_{\Hspace{D}_z}+P_{\Hspace{D}_{*,z}}))^n
P_{\Hspace{D}_z}[\begin{array}{ll}
(z\widetilde{\mathbf{A}})|\Hspace{X} & z\widetilde{\mathbf{B}}
\end{array}]
\end{eqnarray*}
since by \eqref{eq:n-inv'}
$P_\Hspace{X}(z\widetilde{\mathbf{A}})P_{\Hspace{D}_z}=0,\
(z\widetilde{\mathbf{C}})P_{\Hspace{D}_z}=0,\
P_{\Hspace{D}_{*,z}}(z\widetilde{\mathbf{A}})|\Hspace{X}=0$ and
$P_{\Hspace{D}_{*,z}}(z\widetilde{\mathbf{B}})=0$. Further, by
\eqref{eq:n-inv'}
$z\widetilde{\mathbf{A}}\Hspace{D}_z\subset\Hspace{D}_z\perp
\Hspace{D}_{*,z}$, hence
\begin{displaymath}
P_{\Hspace{D}_{*,z}}
(z\widetilde{\mathbf{A}}(P_{\Hspace{D}_z}+P_{\Hspace{D}_{*,z}}))^n
P_{\Hspace{D}_z}=P_{\Hspace{D}_{*,z}}
(z\widetilde{\mathbf{A}})^nP_{\Hspace{D}_z}=0,
\end{displaymath}
and therefore $z\mathbf{H}(z\mathbf{T})^nz\mathbf{F}=0$, i.e. the
necessary condition of this theorem is fulfilled.

Conversely, let $L_\mathbf{G}$ belong to the class
$S_N^0(\Hspace{X}\oplus
\Hspace{N^-},\Hspace{X}\oplus\Hspace{N^+})$. Then by
Theorem~\ref{thm:tf-cons} there exists such a conservative LSDS
$\beta=(N; \mathbf{T}, \mathbf{F}, \mathbf{H}, \mathbf{S};
\Hspace{Y}, \Hspace{X}\oplus\Hspace{N^-},
\Hspace{X}\oplus\Hspace{N^+})$ that for all $z\in\nspace{D}{N}$
\begin{displaymath}
\theta_\beta
(z)=z\mathbf{S}+z\mathbf{H}(I_\Hspace{Y}-
z\mathbf{T})^{-1}z\mathbf{F}=z\mathbf{G}=L_\mathbf{G}(z).
\end{displaymath}
Then $\mathbf{S}=\mathbf{G}$, and for all $n\in\mathbb{Z_+}$ and
$z\in\nspace{D}{N}$ (and hence for all $z\in\nspace{C}{N}$)
\begin{equation}\label{eq:vanish}
z\mathbf{H}(z\mathbf{T})^nz\mathbf{F}=0.
\end{equation}
Conservativity of $\beta$ means that for each $\zeta\in\nspace{T}{N}$
\begin{displaymath}
\zeta\widetilde{\mathbf{G}}=\left(
\begin{array}{ll}
\zeta\mathbf{T} & \zeta\mathbf{F} \\
\zeta\mathbf{H} & \zeta\mathbf{S}
\end{array}
\right)=\left(
\begin{array}{ll}
\zeta\mathbf{T} & \zeta\mathbf{F} \\
\zeta\mathbf{H} & \zeta\mathbf{G}
\end{array}
\right)
\in [\Hspace{Y}\oplus
(\Hspace{X}\oplus
\Hspace{N^-}),\Hspace{Y}\oplus
(\Hspace{X}\oplus\Hspace{N^+})]
\end{displaymath}
is a unitary operator. This operator allows another block partition:
\begin{displaymath}
\zeta\widetilde{\mathbf{G}}:=\left(
\begin{array}{ll}
\zeta\widetilde{\mathbf{A}} & \zeta\widetilde{\mathbf{B}} \\
\zeta\widetilde{\mathbf{C}} & \zeta\widetilde{\mathbf{D}}
\end{array}
\right)\in [(\Hspace{Y}\oplus\Hspace{X})\oplus
\Hspace{N^-},(\Hspace{Y}\oplus\Hspace{X})\oplus\Hspace{N^+}]
\end{displaymath}
where
\begin{equation}\label{eq:cons-n-dil}
\begin{array}{ll}
\zeta\widetilde{\mathbf{A}}=\left[
\begin{array}{cc}
\zeta\mathbf{T} & (\zeta\mathbf{F})|\Hspace{X} \\
P_\Hspace{X}(\zeta\mathbf{H}) & \zeta\mathbf{A}
\end{array}
\right], &
\zeta\widetilde{\mathbf{B}}=\left[
\begin{array}{c}
(\zeta\mathbf{F})|\Hspace{N^-} \\
\zeta\mathbf{B}
\end{array}
\right], \\
\zeta\widetilde{\mathbf{C}}=
[\begin{array}{cc}
P_\Hspace{N^+}(\zeta\mathbf{H}) & \zeta\mathbf{C}
\end{array}], & \zeta\widetilde{\mathbf{D}}=\zeta\mathbf{D}.
\end{array}
\end{equation}
It is clear that one can correspond to this partition of
$\zeta\widetilde{\mathbf{G}}$ the conservative scattering LSDS
$\widetilde{\alpha}
=(N; \widetilde{\mathbf{A}}, \widetilde{\mathbf{B}},
\widetilde{\mathbf{C}}, \mathbf{D};\Hspace{Y}\oplus\Hspace{X},
\Hspace{N^-},
\Hspace{N^+})$. Let us show that $\widetilde{\alpha }$ is a dilation of
$\alpha $. For this purpose, according to 
Proposition~\ref{prop:reform-n-dil}, it is sufficient to verify the equalities
in \eqref{eq:n-mom} 
for all $\zeta\in\nspace{T}{N}$. According to \eqref{eq:cons-n-dil}
$\zeta\mathbf{A}=P_\Hspace{X}(\zeta\widetilde{\mathbf{A}})|\Hspace{X}$,
i.e. for $n=1$ the first equality in \eqref{eq:n-mom} holds (for $n=0$
it holds trivially). Let us apply induction on $n$. Suppose that
$(\zeta\mathbf{A})^n=P_\Hspace{X}(\zeta\widetilde{\mathbf{A}})^n|
\Hspace{X}$ for $n=k\in\mathbb{Z_+}\setminus\{ 0\}$. Then by
\eqref{eq:cons-n-dil} and \eqref{eq:vanish} we have
\begin{eqnarray*}
(\zeta\mathbf{A})^{k+1} &=&
(\zeta\mathbf{A})(\zeta\mathbf{A})^k=P_\Hspace{X}(\zeta
\widetilde{\mathbf{A}})P_\Hspace{X}
(\zeta\widetilde{\mathbf{A}})^k|\Hspace{X} \\
&=&
P_\Hspace{X}(\zeta\widetilde{\mathbf{A}})(I_{\Hspace{Y}\oplus\Hspace{X}}
-P_\Hspace{Y})(\zeta\widetilde{\mathbf{A}})^k|\Hspace{X} \\
&=&
P_\Hspace{X}(\zeta\widetilde{\mathbf{A}})^{k+1}|\Hspace{X}-
P_\Hspace{X}(\zeta\mathbf{H})P_\Hspace{Y}(\zeta\widetilde{\mathbf{A}})^k
|\Hspace{X} \\
&=& P_\Hspace{X}(\zeta\widetilde{\mathbf{A}})^{k+1}|\Hspace{X}-
P_\Hspace{X}(\zeta\mathbf{H})
[\begin{array}{cc}
\zeta\mathbf{T} & (\zeta\mathbf{F})|\Hspace{X}
\end{array}](\zeta\widetilde{\mathbf{A}})^{k-1}|\Hspace{X} \\
&=& P_\Hspace{X}(\zeta\widetilde{\mathbf{A}})^{k+1}|\Hspace{X}-
P_\Hspace{X}(\zeta\mathbf{H})(\zeta\mathbf{T})P_\Hspace{Y}(\zeta
\widetilde{\mathbf{A}})^{k-1}|\Hspace{X}=\ldots \\
&=& P_\Hspace{X}(\zeta\widetilde{\mathbf{A}})^{k+1}|\Hspace{X}-
P_\Hspace{X}(\zeta\mathbf{H})(\zeta\mathbf{T})^{k-
1}P_\Hspace{Y}(\zeta\widetilde{\mathbf{A}})|\Hspace{X} \\
&=& P_\Hspace{X}(\zeta\widetilde{\mathbf{A}})^{k+1}|\Hspace{X}-
P_\Hspace{X}(\zeta\mathbf{H})(\zeta\mathbf{T})^{k-
1}(\zeta\mathbf{F})|\Hspace{X}=P_\Hspace{X}(\zeta\widetilde{\mathbf{A}})
^{k+1}|\Hspace{X}.
\end{eqnarray*}
Thus the first equality in \eqref{eq:n-mom} is valid for all
$n\in\mathbb{Z_+}$. The second and the third equalities in
\eqref{eq:n-mom} are proved analogously. Finally, for an arbitrary
$n\in\mathbb{Z_+}$ by \eqref{eq:cons-n-dil}, \eqref{eq:vanish} and the
second equality in \eqref{eq:n-mom} we have
\begin{eqnarray*}
\zeta\mathbf{C}(\zeta\mathbf{A})^n\zeta\mathbf{B} &=&
(\zeta\widetilde{\mathbf{C}})P_\Hspace{X}(\zeta\widetilde{\mathbf{A}})^n
\zeta\widetilde{\mathbf{B}}=(\zeta\widetilde{\mathbf{C}})(I_{\Hspace{Y}
\oplus\Hspace{X}}-P_\Hspace{Y})(\zeta\widetilde{\mathbf{A}})^n\zeta
\widetilde{\mathbf{B}} \\
&=&
\zeta\widetilde{\mathbf{C}}(\zeta\widetilde{\mathbf{A}})^n
\zeta\widetilde{\mathbf{B}}-P_\Hspace{N^+}(\zeta\mathbf{H})P_\Hspace{Y}
(\zeta\widetilde{\mathbf{A}})^n\zeta\widetilde{\mathbf{B}}\\
&=&
\zeta\widetilde{\mathbf{C}}(\zeta\widetilde{\mathbf{A}})^n
\zeta\widetilde{\mathbf{B}}-
P_\Hspace{N^+}(\zeta\mathbf{H})[\begin{array}{cc}
\zeta\mathbf{T} & (\zeta\mathbf{F})|\Hspace{X}
\end{array}](\zeta\widetilde{\mathbf{A}})^{n-
1}\zeta\widetilde{\mathbf{B}}\\
&=& \zeta\widetilde{\mathbf{C}}(\zeta\widetilde{\mathbf{A}})^n
\zeta\widetilde{\mathbf{B}}-
P_\Hspace{N^+}(\zeta\mathbf{H})(\zeta\mathbf{T})P_\Hspace{Y}
(\zeta\widetilde{\mathbf{A}})^{n-1}\zeta\widetilde{\mathbf{B}}=\ldots \\
&=& \zeta\widetilde{\mathbf{C}}(\zeta\widetilde{\mathbf{A}})^n
\zeta\widetilde{\mathbf{B}}-
P_\Hspace{N^+}(\zeta\mathbf{H})(\zeta\mathbf{T})^nP_\Hspace{Y}
\zeta\widetilde{\mathbf{B}}\\
&=& \zeta\widetilde{\mathbf{C}}(\zeta\widetilde{\mathbf{A}})^n
\zeta\widetilde{\mathbf{B}}-
P_\Hspace{N^+}(\zeta\mathbf{H})(\zeta\mathbf{T})^n(\zeta\mathbf{F})|
\Hspace{N^-}=\zeta\widetilde{\mathbf{C}}(\zeta\widetilde{\mathbf{A}})^n
\zeta\widetilde{\mathbf{B}}.
\end{eqnarray*}
Note that for $n=0$ this calculation is obviously simplified and does
not contain terms like $(\zeta\mathbf{A})^k$ with $k<n$. The proof is
complete.
\end{proof}
In the particular case when $\Hspace{N^-}=\Hspace{N^+}=\{ 0\}$ we obtain
the following result.
\begin{cor}
The linear pencil of contractions $L_\mathbf{A}(\zeta
):=\zeta\mathbf{A}\in [\Hspace{X}, \Hspace{X}]\ (\zeta\in\nspace{T}{N})$
allows a \emph{unitary dilation}, i.e. there is a linear pencil of unitary
operators $L_{\widetilde{\mathbf{A}}}(\zeta
)=\zeta\widetilde{\mathbf{A}}\in [\Hspace{X}, \Hspace{X}]\
(\zeta\in\nspace{T}{N})$, for which
$\Hspace{X}\subset\widetilde{\Hspace{X}}$ and
\begin{equation}\label{eq:pencil-dil}
\forall\zeta\in\nspace{T}{N},\ \forall n\in\mathbb{Z_+}\quad
(\zeta\mathbf{A})^n=P_\Hspace{X}(\zeta\widetilde{\mathbf{A}})^n|\Hspace{
X},
\end{equation}
if and only if $L_\mathbf{A}\in S_N^0(\Hspace{X},\Hspace{X})$.
\end{cor}
\begin{rem}
Condition \eqref{eq:pencil-dil} is equivalent to the family of equalities in
the first line of condition
\eqref{eq:sym-mom} in Proposition~\ref{prop:power-n-dil}.
\end{rem}

\section{$N$-dissipative scattering LSDSs}\label{sec:n-pass}
It is obvious (see Section~\ref{sec:pr}) that the multiparametric LSDS
$\alpha =(N; \mathbf{A}, \mathbf{B}, \mathbf{C},
\mathbf{D};\Hspace{X},
\Hspace{N^-}, \Hspace{N^+})$ is a dissipative scattering LSDS if and
only if
the corresponding linear function $L_\mathbf{G}$ in
\eqref{eq:linear-func} belongs to the class $B_N^0(\Hspace{X}\oplus
\Hspace{N^-},\Hspace{X}\oplus\Hspace{N^+})$.
\begin{defn}\label{defn:n-pass}
We shall call the system $\alpha =(N; \mathbf{A}, \mathbf{B},
\mathbf{C}, \mathbf{D};\Hspace{X}, \Hspace{N^-}, \Hspace{N^+})$ a
$N$-\emph{dissipative scattering LSDS} if $L_\mathbf{G}\in
S_N^0(\Hspace{X}\oplus
\Hspace{N^-},\Hspace{X}\oplus\Hspace{N^+})$.
\end{defn}
It is clear that by virtue of \eqref{eq:inclus} the class of
$N$-dissipative scattering LSDSs is a subclass of the class of all
dissipative
scattering LSDSs. According to Theorem~\ref{thm:criterion} it consists
of those and only those dissipative systems which allow conservative
dilations. It follows from Theorem~\ref{thm:tf-pass} that the class of
transfer functions of $N$-parametric dissipative scattering LSDSs with
the
input space $\Hspace{N^-}$ and the output space $\Hspace{N^+}$ is a
subclass of $B_N^0(\Hspace{N^-},\Hspace{N^+})$, however we have no
complete description of this subclass. For $N$-dissipative systems, from
Theorem~\ref{thm:tf-cons}, Theorem~\ref{thm:criterion} and
Proposition~\ref{prop:tf-n-dil} we obtain the complete description of
the class of transfer functions.
\begin{thm}
The class of transfer functions of $N$-dissipative scattering LSDSs with
the input space $\Hspace{N^-}$ and the output space $\Hspace{N^+}$
coincides with $S_N^0(\Hspace{N^-}, \Hspace{N^+})$.
\end{thm}
Let us note that, by Theorem~\ref{thm:criterion}, if the $N$-dissipative
scattering LSDS $\alpha $ is a dilation of some system $\alpha _0$
then $\alpha _0$ is also $N$-dissipative. By virtue of
Theorem~\ref{thm:tf-cons}, for each operator-valued function
$\theta\in S_N^0(\Hspace{N^-}, \Hspace{N^+})$ there exists a
\emph{conservative realization} i.e. such a conservative scattering LSDS
$\alpha $ that $\theta =\theta_\alpha $. According to
Proposition~\ref{prop:exist-min}, for $\alpha $ there exists a minimal
system $\alpha _{min}$ such that $\alpha $ is a dilation of $\alpha
_{min}$, moreover by Proposition~\ref{prop:tf-n-dil} $\theta
=\theta_\alpha =\theta_{\alpha _{min}}$. Thus we obtain the theorem on a
\emph{minimal} $N$-\emph{dissipative realization} for operator-valued
functions from $S_N^0(\Hspace{N^-}, \Hspace{N^+})$.
\begin{thm}
For an arbitrary $\theta\in S_N^0(\Hspace{N^-}, \Hspace{N^+})$ there
exists a minimal $N$-dissipative scattering LSDS $\alpha_{min}$ such
that
$\theta =\theta_{\alpha _{min}}$.
\end{thm}
As we remarked in the end of Section~\ref{sec:pr}, for $N=1$ and $N=2$
we have equality in \eqref{eq:inclus}. It follows from here (see
Definition~\ref{defn:n-pass}) that for these cases the notions of
$N$-dissipative scattering LSDS and dissipative scattering LSDS
coincide, and
the corresponding classes of systems also coincide. In the case $N>2$,
as we will show, these classes do not coincide.

The following result was obtained in \cite{K1}.
\begin{thm}\label{thm:anti-vN}
There exist such triples $\mathbf{T}=(T_1, T_2, T_3)$ of commuting
contractions on some finite-dimensional Hilbert space $\Hspace{H}$ and
$\mathbf{M}=(M_1, M_2, M_3)$ of linear operators on
$\nspace{C}{n}$, with some integer $n>1$, that the linear
homogeneous operator-valued function $L_\mathbf{M}(z_1, z_2,
z_3)=M_1z_1+M_2z_2+M_3z_3$ satisfies
\begin{equation}\label{eq:anti-vN}
\| L_\mathbf{M}(\mathbf{T})\|>\max_{z\in\overline{\nspace{D}{3}}}\|
L_\mathbf{M}(z)\|
\end{equation}
(here $L_\mathbf{M}(\mathbf{T}):=M_1\otimes
T_1+M_2\otimes T_2+M_3\otimes T_3\in [\nspace{C}{n}\otimes \Hspace{H},
\nspace{C}{n}\otimes \Hspace{H}]$).
\end{thm}
Under conditions of this theorem, set
\begin{equation}\label{eq:ex}
G_k:=(\max_{z\in\overline{\nspace{D}{3}}}\| L_\mathbf{M}(z)\| )^{-
1}M_k,\quad (k=1, 2, 3)
\end{equation}
\begin{equation}\label{eq:ex-sys1}
\Hspace{X}:=\nspace{C}{{n-1}},\ \Hspace{N^-}=\Hspace{N^+}:=\mathbb{C},\
\mbox{so that}\ \nspace{C}{n}=\Hspace{X}\oplus\Hspace{N^-
}=\Hspace{X}\oplus\Hspace{N^+},
\end{equation}
\begin{equation}\label{eq:ex-sys2}
\begin{array}{ll}
A_k:=P_\Hspace{X}G_k|\Hspace{X}, & B_k:=P_\Hspace{X}G_k|\Hspace{N^-} \\
C_k:=P_\Hspace{N^+}G_k|\Hspace{X}, & D_k:=P_\Hspace{N^+}G_k|\Hspace{N^-
}.
\end{array}
\quad (k=1, 2, 3)
\end{equation}
Then the linear operator-valued function
$L_\mathbf{G}(z):=G_1z_1+G_2z_2+G_3z_3=z\mathbf{G}\
(z\in\nspace{D}{3})$, by virtue of \eqref{eq:ex}, belongs to the class
$B_3^0(\nspace{C}{n}, \nspace{C}{n})$. However by \eqref{eq:anti-vN} and
\eqref{eq:ex}
\begin{displaymath}
\| L_\mathbf{G}(\mathbf{T})\| > \max_{z\in\overline{\nspace{D}{3}}}\|
L_\mathbf{G}(z)\| =1,
\end{displaymath}
and hence there is a positive $r<1$ for which
\begin{displaymath}
\| L_\mathbf{G}(r\mathbf{T})\| > 1.
\end{displaymath}
The latter means (see Section~\ref{sec:pr}) that $L_\mathbf{G}$ does not
belong to the class $S_3^0(\nspace{C}{n}, \nspace{C}{n})$. Thus the LSDS
$\alpha =(3; \mathbf{A}, \mathbf{B}, \mathbf{C}, \mathbf{D};\Hspace{X},
\Hspace{N^-}, \Hspace{N^+})$, that is defined in
\eqref{eq:ex}--\eqref{eq:ex-sys2}, is dissipative but not 3-dissipative.
For
the case $N>3$ an analogous example of dissipative but not $N$-
dissipative
system can be easily constructed by supplement of arbitrary operators
$M_4,\ldots ,M_N$ on $\nspace{C}{n}$ with sufficiently small norms to
the triple
$\mathbf{M}=(M_1, M_2, M_3)$ from Theorem~\ref{thm:anti-vN} and setting
$T_4=\ldots =T_N=0$, so that the inequality analogous to
\eqref{eq:anti-vN} holds for the $N$-tuples
$\widetilde{\mathbf{M}}:=(M_1, M_2, M_3, M_4,\ldots , M_N)$ and
$\widetilde{\mathbf{T}}:=(T_1, T_2, T_3, 0,\ldots ,0)$,
and then defining such a $N$-parametric system in the same way as in
\eqref{eq:ex}--\eqref{eq:ex-sys2}. Thus we have proved the following.
\begin{thm}
The class of $N$-dissipative scattering LSDSs for the cases $N=1$ and
$N=2$
coincides with the class of all $N$-parametric dissipative scattering
LSDSs,
and for the case $N>2$ is a proper subclass of the latter.
\end{thm}

\bibliographystyle{amsplain}
\bibliography{msys}
\ \\
Department of Higher Mathematics \\
Odessa State Academy of Civil Engineering and Architecture \\
Didrihson str. 4, Odessa, 270029 \\
Ukraine \\
\\
1991 Mathematics Subject Classification: 47A20, 47A56, 93C35

\end{document}